\documentclass[11pt]{article}
in \oddsidemargin 0.25true in \evensidemargin 0.25true in
\usepackage{epsfig}

\usepackage{enumerate}

\begin{document}
\newcommand{\arcsech}{\mbox{arcsech}}
\newcommand{\bc}{\begin{center}}
\newcommand{\ec}{\end{center}}
\newcommand{\be}{\begin{equation}}
\newcommand{\ee}{\end{equation}}
\newcommand{\bea}{\begin{eqnarray}}
\newcommand{\eea}{\end{eqnarray}}
\newcommand{\bean}{\begin{eqnarray*}}
\newcommand{\eean}{\end{eqnarray*}}
\newcommand{\bd}{\begin{description}}
\newcommand{\ed}{\end{description}}

\newcommand{\caln}{\cal N}
\newcommand{\caloi}{{\cal OI}}
\newcommand{\hatv}{{\bf v}}
\newcommand{\hatrho}{\hat{{\bf \rho}}}
\newcommand{\hatlam}{\lambda}
\newcommand{\tillam}{\tilde{\lambda}}
\newcommand{\tT}{\tilde{T}}
\newcommand{\hT}{\hat{T}}
\newcommand{\bT}{{\bf T}}
\newcommand{\bA}{{\bf A}}
\newcommand{\ep}{\epsilon}
\newcommand{\cd}{{\cal D}}
\newcommand{\hp}{\hat{T}}
\newcommand{\epf}{\epsilon_f}
\newcommand{\Nn}{N_{net}}
\newcommand{\Ns}{N_{sph}}
\newcommand{\Nd}{N_{dat}}
\newcommand{\lng}{\langle}
\newcommand{\rng}{\rangle}
\newcommand{\vep}{\epsilon}
\newcommand{\vph}{\varphi}
\newcommand{\vp}{\varphi}
\newcommand{\bfk}{{\bf k}}
\newcommand{\bfr}{{\bf r}}
\newcommand{\bfN}{{\bf N}}
\newcommand{\bfa}{{\bf a}}
\newcommand{\bfb}{{\bf b}}
\newcommand{\bfc}{{\bf c}}
\newcommand{\bfH}{{\bf H}}
\newcommand{\bfK}{{\bf K}}
\newcommand{\bfL}{{\bf L}}
\newcommand{\bfM}{{\bf M}}
\newcommand{\bfx}{{\bf x}}
\newcommand{\bC}{{\bf C}}
\newcommand{\bp}{{\bf{p}}}
\newcommand{\bfe}{{\bf e}}
\newcommand{\hbt}{\hat{\bar{T}}}
\newcommand{\uvr}{{\bf \hat{r}}}
\newcommand{\uvn}{{\bf \hat{n}}}
\newcommand{\bw}{\bar{w}}
\newcommand{\bx}{\mbox{\boldmath $\xi$}}
\newcommand{\bnu}{\mbox{\boldmath $\nu$}}
\newcommand{\my}{M_{yr}}

\newcommand{\bvt}{\begin{verbatim}}
\newcommand{\evt}{\end{verbatim}}
\newcommand{\dbar}{\vrule width 0.25truein height 1.0pt depth -0.6pt}

\renewcommand{\baselinestretch}{2.0}

\vskip 0.5cm

\renewcommand{\baselinestretch}{1.2}

\begin{center}
{\bf\Huge Calculus from a \\
\vskip 0.3cm

Statistics Perspective}
\vskip 0.5cm
\end{center}

\vskip 1.5cm
\begin{center}
{\Large SAMUEL S.P. SHEN, DOV ZAZKIS, KIMBERLY LEUNG, and \\
\vskip 0.3cm

CHRIS RASMUSSEN }  \\
\vskip 0.5cm

Department of Mathematics and Statistics, 
San Diego State University,\\
San Diego, CA 92182, USA. \\Email: 
{\tt sam.shen@sdsu.edu}\\

\vskip 2.0cm
 
\end{center}

\thispagestyle{empty}
\renewcommand{\baselinestretch}{1.2}

\noindent {{\bf  Abstract}} This paper provides an approach to establishing the calculus method from the concept of mean, 
i.e., average. This approach is from a statistics perspective and  can help calculus learners understand  calculus ideas and  analyze a 
function defined by data or sampling values from a given function, rather than 
an explicit mathematical formula. The basics of this approach are two averages: arithmetic mean and graphic mean. The 
arithmetic mean is  used to define integral. Area is used to interpret the meaning 
of an integral. Antiderivative is introduced from integral,
and derivative-antiderivative pair is introduced as a mathematical operation entity. The graphic mean is an average 
speed in an interval and is used to interpret the meaning of a derivative.

\newpage

\section{Introduction} 

\noindent  How can one use elementary statistics method to help explain ideas of calculus? The purpose of this paper is to provide an approach  to developing calculus 
methods from a statistics perspective. 
Our approach is motivated in part by mathematical modeling (ref. [{\bf 5}]), and can be helpful for enhancing the  link between  elementary statistics courses and calculus 
courses. The link  is developed from statistics to calculus and is different from the traditional connection from calculus to 
statistics. 

 A mathematical function  $y=f(x)$ can be represented in 
four ways (see  p11 of  Stewart (2008) [{\bf 10}]): (i) a description in words or text, (ii) a table of at least two columns of data, (iii) a graph in 
the x-y coordinates plane, and (iv) an explicit mathematical formula. Conventional calculus textbooks (e.g., [{\bf 1, 10}])  
almost exclusively deal with Case (iv):  functions defined by explicit formulas. The graphic representations (i.e., Case (iii)) are dealt with 
as a supporting tool for the function formula. Case (ii), data-represented functions, is usually considered 
the category of statistics, which analyzes data and makes inferences for the data's implications. Hence, most calculus textbooks 
do not deal with this case even in the modern  era of computers.  However, the rapid development of data-based 
information in our practical life requires us to broaden the traditional statistics methods, 
 including interconnections between statistics and calculus. Today's  speed and availability of
 personal computers and smart phones make
it possible to take new and innovative approaches to calculus for functions  represented in any of the four ways. This paper will develop the calculus ideas from statistics perspective. We will build the calculus concept for the function of Case (ii) and make the concept applicable to the functions of all 
 the four cases.  We use two 
types of averages: arithmetic mean and graphic mean.  The arithmetic mean and law of large numbers (LLN)  are used to define an integral, and graphic mean is to 
define a derivative. 

We do not intend to use this approach to replace the conventional way of teaching 
calculus. Instead, our proposed approach may be used as  supplemental material in today's classrooms of 
conventional calculus and hence provides a new perspective for explaining calculus ideas to students. Further, 
we do not claim that this approach is superior to other approaches to calculus. Every approach to  calculus has its own disadvantages. 
Ours is not an exception. 

Two co-authors have experimented the new approach in two different courses. Shen used  it in  Calculus I for engineering and physical 
sciences in the Fall 2011, where he spent two hours altogether at the beginning of the course on the new approach in a class of 120 students. After the two hours, he asked his students to explain the concept of integral to their grandmother. Zakis 
taught this material to four students  in the summer of 2011 in the third mathematics course in a sequence of four courses for prospective elementary school teachers. This third mathematics course's main contents are on probability and statistics. He made an educational pedagogy study from this course and 
spent over four hours on this approach with emphasis on understanding the  concept of integrals (Zazkis et al. (2014) [ref. {\bf 12}]). 
Therefore, our modest goals of the calculus from statistics perspective are (a) to provide a 
set of supplemental materials  to explain the the concept of calculus from a non-traditional perspective, and (b) to provide research opportunities for mathematics education. 

To simplify the description of our approach to the integral 
and derivative concept, we limit our functions to those that are continuous and smooth with positive domain and positive range, when discussing
the Case (iv) function $y=f(x)$. These limits 
do not affect the rigor of the mathematics developed here and can be easily removed   when 
a more sophisticated mathematics of calculus is introduced.

\section{Arithmetic mean, sampling, and average of a function}

\noindent  For a given location on Earth, its temperature's variation with respect to time forms a functional relationship. 
Consider the arithmetic mean for  the annual surface air temperature (SAT) data (in units  [$^\circ$C])  data from 1951-2010 at Fredericksburg (38.32$^\circ$N, 77.45$^\circ$W), which is 80 kilometers southwest of 
Washington DC, United States (See Table 1).  The data pairs $(x_i, f_i) (i=1, 2, \cdots, n)$ represent
the tabular form of a function (the Case (ii) function) with $x_i$ for time and $f_i$ for temperature, and $n$ is 60 for this case. {\it Arithmetic mean} is defined as 
\be
\bar{f} = \frac{\sum_{i=1}^n f_i}{n}. 
\ee
For the SAT data in Table 1, the mean is 13.2[$^\circ$C].

 \begin{center} 
 Table 1.  Annual mean SAT at Fredericksburg from 1951- 2010
\end{center}
\vskip 0.0in

 \begin{tabular}{| ll| ll| ll| ll| ll| ll|}
 \hline
 \multicolumn{12}{|l|}{Annual mean surface air temperature at Fredericksburg (38.32$^\circ$N, 77.45$^\circ$W)} \\
  \multicolumn{12}{|l|}{Virginia, United States. The temperature units are [$^\circ$C]. The data are from }\\
     \multicolumn{12}{|l|}{the U.S. Historical Climatology Network.}\\
    \multicolumn{12}{|l|}{http://www.ncdc.noaa.gov/oa/climate/research/ushcn/}\\

\hline
 1951&13.1&1961&12.9&1971&13.2&1981&12.5&1991&14.3 &2001&13.7\\
 1952&13.4&1962&12.3&1972&12.8&1982&12.7&1992&12.5&2002&14.3\\
 1953&14.4&1963&12.4&1973&13.7&1983&12.8&1993&13.2&2003&13.0\\
 1954&13.7&1964&13.2&1974&13.3&1984&13.0&1994&13.2&2004&13.6\\
 1955&12.6&1965&12.7&1975&13.0&1985&13.6&1995&13.2&2005&13.7\\
 1956&13.1&1966&12.4&1976&12.6&1986&13.4&1996&12.6&2006&14.2\\
 1957&13.3&1967&12.3&1977&13.1&1987&13.4&1997&13.0&2007&14.1\\
 1958&12.1&1968&12.8&1978&12.3&1988&12.6&1998&14.7&2008&13.7\\
 1959&13.6&1969&12.5&1979&12.5&1989&12.9&1999&13.9&2009&13.2\\
 1960&12.6 &1970&13.1&1980&12.7&1990&14.3&2000&13.0&2010&14.2\\
 \hline
 \end{tabular}

 \vskip 0.5cm

Next we explore the data that are samples from the Case (iv) functional values. 
The samples of independent and dependent variables are $(x_i, y_i), i=1.2.\cdots, n$, where the Case (iv) 
function is denoted by $y=f(x)$.  The sampling of $x$ is
usually done in three ways.
\begin{enumerate}[(i)]
\item Uniform sampling: The distance between each pair of neighborhood samples is the same, i.e., $h=x_{i}-x_{i-1}$ 
is the same for each $i$. The sample points can be determined by $x_{i}=x_{i-1} + h, i=1, 2, \cdots, n$. 
An example is the SAT sampling in Table 1 whose $h$ is one year. Another example is a uniform sampling of a function shown 
in Figure 1. 
\item Random sampling: Random sampling is, by name, a probability process and can be determined by  a predefined way using a statistics 
software. Many software tools for  generating random numbers are available in pubic domain for free, such as WolframAlpha:
{\tt www.wolframalpha.com}.
The command {\tt RandomReal[\{0,2\},10] }yields 10 random real numbers in the interval $[0,2]$. A trial of this command is\\
  \small{
  \{1.25685, 1.02584, 1.26752, 0.813785, 0.224482, 0.905491, 1.73265, 1.44104, 1.95108, 0.926807\}
  }\\
  Of course, each trial yields a different random result. 
\item Convenience sampling: This sampling may be neither as evenly spaced as the uniform sampling nor completely random, but 
the samples are taken from convenient and practical  locations, such as the location of weather stations, which could  be neither
at the top of Himalayas nor over the midst of Pacific. Thus, convenient sampling is often encountered in practice, such as 
the data from geology, meteorology, hydrology, agriculture and forestry, but suffers possible drawbacks of bias and incompletion. 
When given a Case (iv) function, this sampling is unnecessary, but for a function of the other three cases in practical applications, 
this sampling is essential to describe a function for applications. 

\end{enumerate}
 \begin{figure} [ht]
\centering
\includegraphics[height=3.5in, width=4.5in]{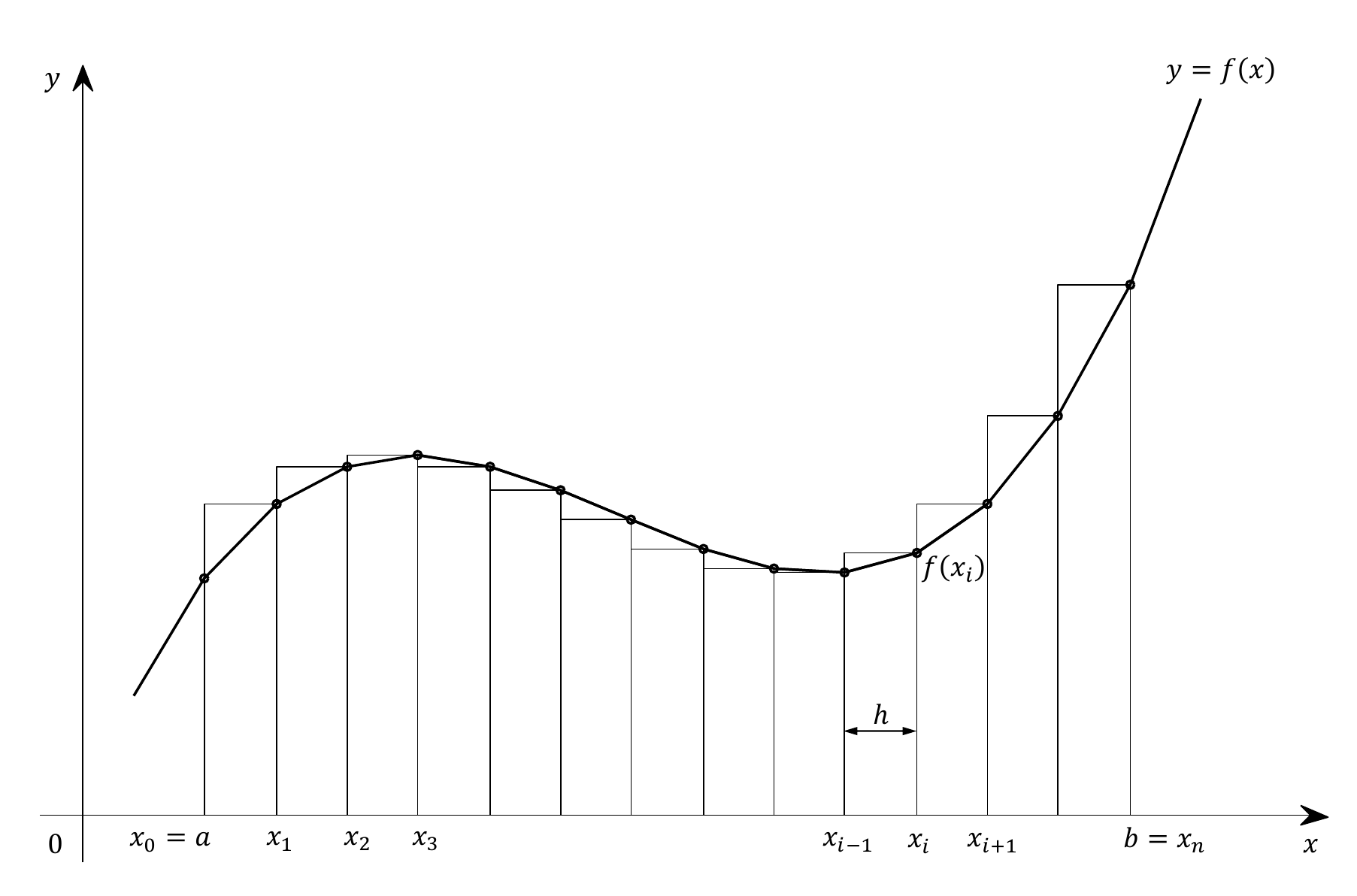}
\caption{Uniform sampling for a function $y=f(x)$ in $[a,b]$. The area under the curve is approximated by the sum of the
narrow rectangular strips}
\label{figure1}
\end{figure}
\vskip 0.3cm

Let us explore two examples of calculating arithmetic means of sample data from a function: $y=x$ and $y=x^2$.

From our intuition or assisted by a figure (see Figure 2 with $a=0$), we can infer that the mean of the function  $y=x$ in $[0,1]$ is $1/2$. 
This result can be simulated by different sampling methods. 
\begin{enumerate}[(a)]
\item Uniform sampling: We use a sample of size 100 with the first sample at $x=0.01$ and last sample at $x=1$. 
These 100 points divide the interval $[0,1]$ into 100 equal sub-intervals of length $1/100$:
\[
x_1=0.01, x_2=0.02, \cdots, x_{i}= \frac{i}{100}, \cdots, x_{100}=1.
\]
\[
 y_i=x_i, i= 1,2,\cdots, 100
 \]
The arithmetic mean of data  $y_i, i=0,1, \cdots, 100$ is equal to 
\be
\hat{\bar{y}}=  \frac{\sum_{i=1}^{100} y_i}{100}=\frac{\sum_{i=1}^{100}(i/100)}{100}=\frac{\sum_{i=1}^{100}i}{100\times100}= 0.505,
\ee
which is very close to the exact average value 0.5. 
\item Random sampling: Random sampling usually does not sample the end points. Our 100 samples are at the 
internal points of $[0,1]$, i.e., over $(0,1)$. WolframAlpha command {\tt RandomReal[\{0,1\},100] } generates 100 random 
numbers whose  mean is $0.5312$, which is 
 close to the exact result $0.5$. Again, the result $0.5312$ is different each time the command is implemented due to the random 
 nature of the {\tt RandomReal} generator. It is intuitive that larger samples should be likely to lead to more accurate approximation. When using 1000 samples, a result is $0.4981$, which is very close to $0.5$.
\end{enumerate}

\begin{figure} [ht]
\centering
\includegraphics[height=3in, width=4in]{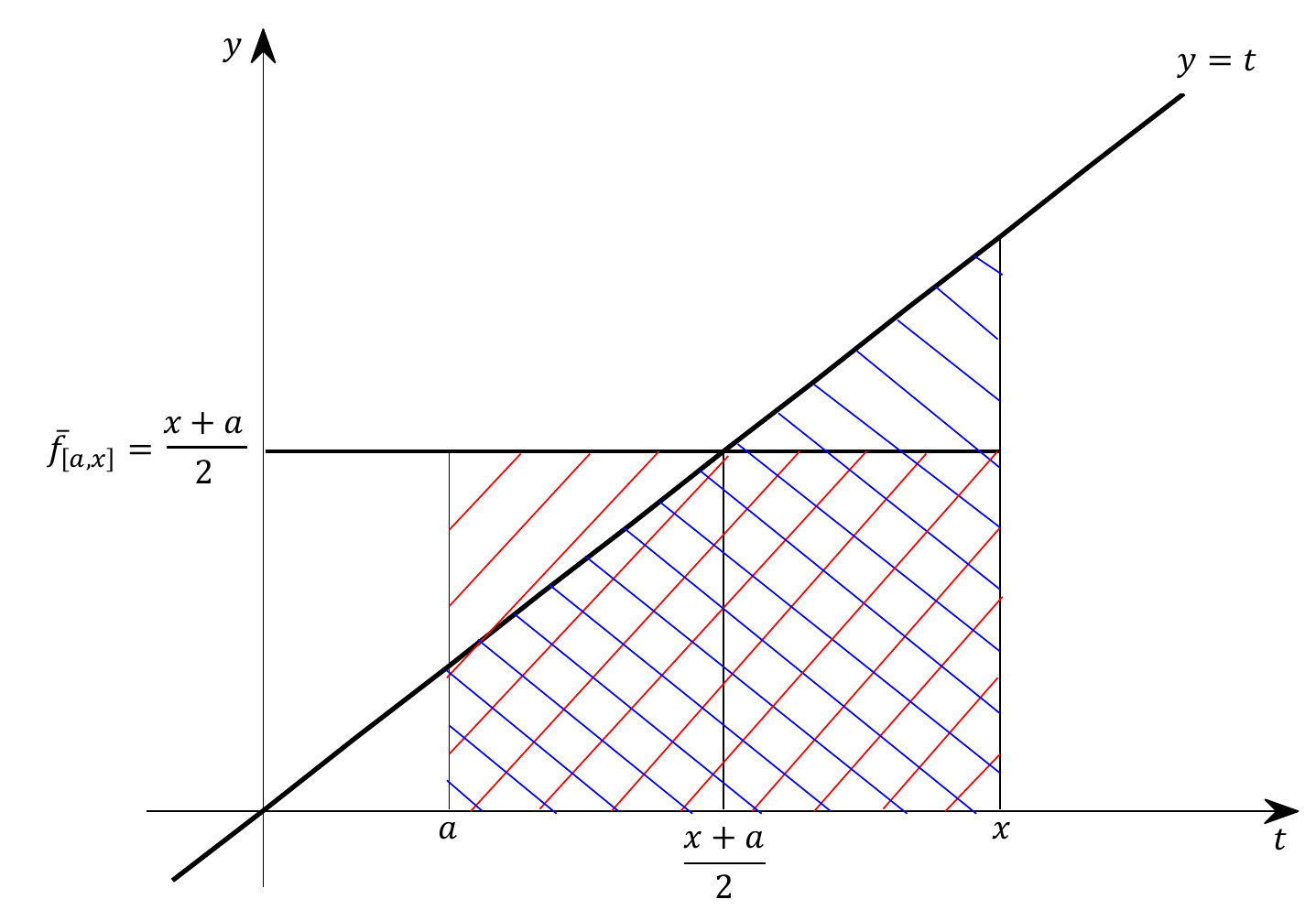}
\caption{ Arithmetic mean of a linear function.}
\label{figure2}
\end{figure}

Next we consider $y=x^2$ in $[0,1]$ (see Figure 3). The uniform sampling of using 100 points as above yields the following:
\[
\hat{\bar{y}}= \frac{\sum_{i=1}^{100}y_i^2}{100}=\frac{\sum_{i=1}^{100}(i/100)^2}{100}=\frac{\sum_{i=0}^{100}i^2}{100^2\times100}= \frac{100\times(100+1)(2\times 100 +1)/6}{100^2\times 100}
= \frac{101\times 201}{6\times 100^2}= .33835
\]
\\
The random sampling with WolframAlpha uses a command 
{\tt RandomReal[\{0,1\},10]}\^{\tt 2}
to generate the data, and the text result can be copied and used to calculate the mean by using WolframAlpha.
A result is $0.36726$. The accurate mean value should be $0.33333$,  which can be computed from $n$ samples when $n$ is very large.

\begin{figure} [ht]
\centering
\includegraphics[height=3in, width=4in]{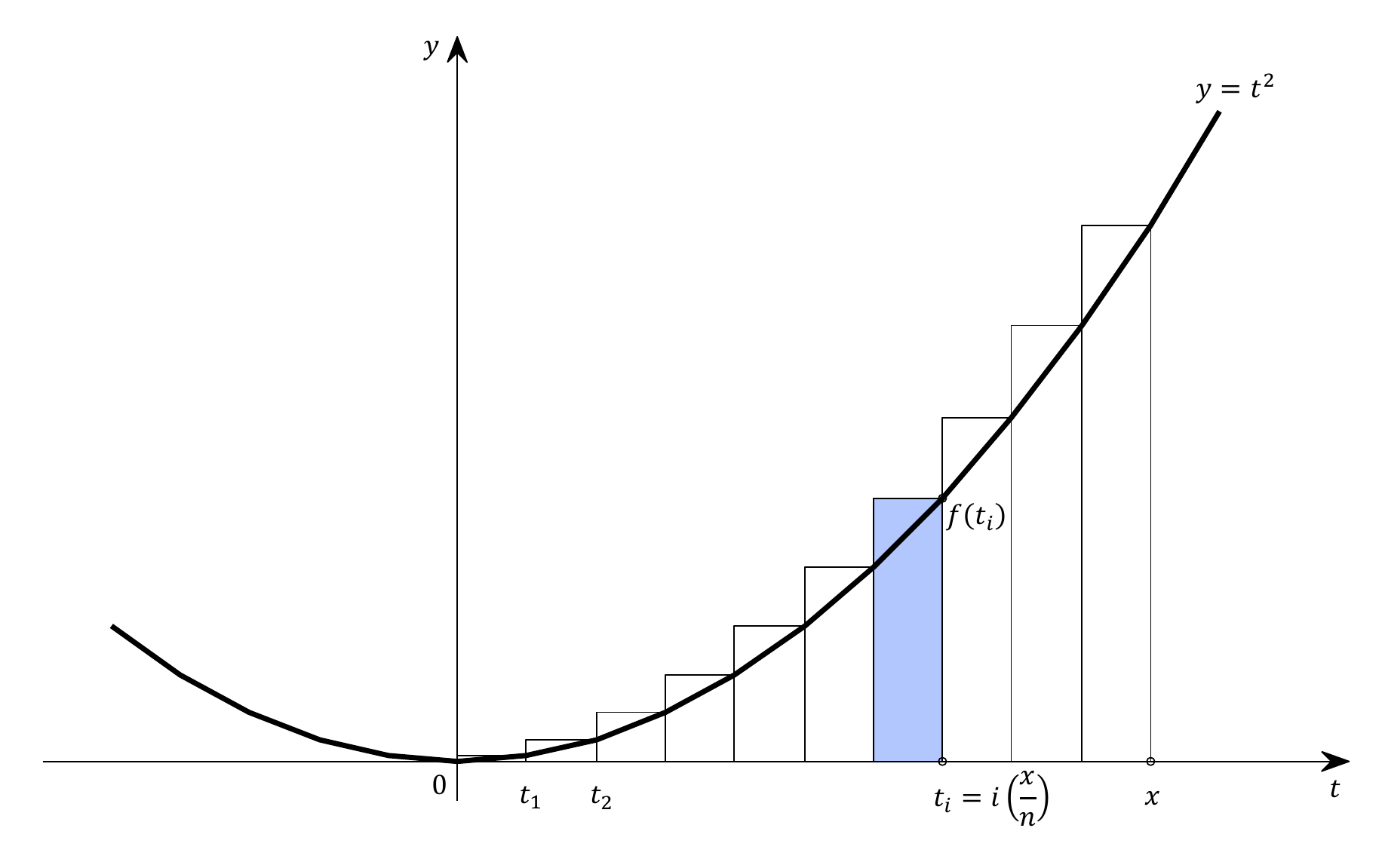}
\caption{Uniform sampling of a parabolic function.}
\label{figure3}
\end{figure}

WolframAlpha has limited statistical computing power. The open source  statistics software R and the commonly used 
MS Excel have more statistical computing power and can be used for practical applications and research. 
An Excel calculation of the sample mean for $y=x^2$ in $[0,1]$ for different number of samples is given in Table 2, which indicates that, in general, the accuracy of the mean improves as the sample size $n$ increases. As a matter of fact, the improvement 
 of the accuracy can be quantified by the error $\sigma/\sqrt{n}$ where $\sigma$ is the 
 standard deviation of the population being sampled. This is a result of the Central Limit Theorem (CLT) in statistics (see a basic statistics text, e.g., [{\bf 4, 11}] ). 
 The CLT states that the mean of samples with sufficiently large sample size is normally distributed. 
 CLT further asserts that (i) the expected value of the sample mean is the same as the population mean,
 and (ii) the variance 
 of the sample mean is $(1/n)$th of the population variance. These two assertions are  widely assumed  in practical
 statistical applications and computer calculations, such as Monte Carlo simulations. 
 This is normal 
 distribution result can be intuitive to students and 
 may be accepted as an axiom. However, instructor does not need to teach the CLT to students in this first introduction of calculus 
 from statistics perspective. This is like the case that the first introduction of conventional calculus does not need to prove the 
 existence of a limit from the rigorous $\delta-\epsilon$ argument. 
 
 \vskip 0.5cm
 
   \begin{center} 
 Table 2. Convergence of the sample mean as the sample size increases. 
\end{center}

 \begin{tabular}{|c |c |c |c |c |c |c |}
 \hline
Sample Size {\it n}&10&100&1,000&10,000&100,000&1,000,000\\
\hline
Uniform Sampling&0.3850&0.3384&0.3338&0.3334&0.3333&0.3333\\
\hline
Random Sampling &\multicolumn{6}{|l|}{}\\
 \hline
 Trial 1&0.4350&0.3505&0.3463&0.3344&0.3315&0.3332\\
 \hline
 Trial 2&0.3560&0.3058&0.3284&0.3325&0.3331&0.3332\\
 \hline
 Trial 3&0.4640&0.3518&0.3355&0.3317&0.3357&0.3332\\
 \hline
 Average&0.4183&0.3360&0.3367&0.3329&0.3335&0.3332\\
 \hline
 \end{tabular}
 \vskip 0.5cm
 
  Let $\hat{\bar{f}}[n]$ denote the mean from $n$ samples. 
 Then, it is almost always true that $\hat{\bar{f}}[n]$ approaches the true mean of the function as $n$ increases. 
 This convergence has a  probability equal to one. The probability for this to be false is zero. This intuitive conclusion is the strong LLN (see [{\bf 10}]), which asserts that the event of the following limit being true has a probability equal to one:
 \be
 \lim_{n\to \infty} \hat{\bar{f}}[n] =  \bar{f}.
 \ee
Here $ \bar{f}$ is the limit of the sequence $\{\hat{\bar{f}}[1], \hat{\bar{f}}[2], \hat{\bar{f}}[3], \cdots\}$, or simply $\{\hat{\bar{f}}[n]\}$.

For example, when $x^2$ is uniformly sampled by $n$ points over $[0,1]$, 
the mean is 
\[
\hat{\bar{f}}[n]= \frac{\sum_{i=1}^{n}(i/n)^2}{n}=\frac{\sum_{i=0}^{n}i^2}{n\times n^2}= \frac{n(n+1)(2 n +1)/6}{n^3}
=\frac{n+1}{n}\times\frac{2n+1}{6n}.
\]
As $n\to \infty$, the first factor goes to 1 and the second to 2/6=1/3. Thus, the above  mean approaches $1/3\approx 0.33333333$ as $n\to \infty$. 

The above leads to the definition of {\it average of a function} $y=f(x)$ in an interval $[a,b]$ 
\be
\bar{f}=\lim_{n\to \infty}\frac{\sum_{i=1}^n f(x_i)}{n},
\ee
where $\{x_i\}_{i=1}^n$ are sampled from $[a,b]$ by uniform sampling, random sampling, and possibly convenience sampling. We also call $\bar{f}$
{\it arithmetic mean}, or just {\it mean},  of the function $y=f(x)$.

\section{ Definition of integral, antiderivative, and DA pair} 

\noindent {\bf Definition of integral:} {\it If $\bar{f}$ is the average of the function $f(x)$ over the interval $[a,b]$, then  $(b-a)\times  \bar{f}$ is defined 
as the  integral of $f(x)$ over $[a,b]$, denoted by 
\be
I[f,a,b]= (b-a) \bar{f}
\ee
}

This serves as the definition of the definite integral in conventional calculus. Since we will not introduce the concept
of indefinite integral, we thus treat "integral" here as the "definite integral."

From the examples of  mean  in the above section, we can calculate the following integrals:
\[
I[x,0,1]=(1-0)\times (1/2) =1/2,
\]
and
\[
I[x^2,0,1]=(1-0)\times (1/3)=1/3.
\]

\vskip 0.3cm

Then, what is the graphic meaning of integral from the above definition of integral? 
We use uniform sampling to make an interpretation. Since $\hat{\bar{f}}[n] \approx \bar{f}$, we have 
\[
I[f,a,b]\approx  (b-a)\hat{\bar{f}}[n] = (b-a)\frac{ \sum_{i=1}^n f_i}{n} = \sum_{i=1}^n\frac{b-a}{n} f_i. 
\]
Here $f_i= f(x_i)$ is the sample value of the function at the sampling location $x_i$ (see Figure 1). Since this is 
uniform sampling, $x_{i+1}-x_i = h = (b-a)/n$. Thus, each term in the above sum
\[
\frac{b-a}{n} f_i
\]
is the area of a rectangular strip with base $h=(b-a)/n$ and height $f_i$, as shown in Figure 1.
Thus, $(b-a)\hat{\bar{f}}[n]$ is the area of under the echelon, formed by $n$ rectangular strips whose heights 
are determined by the function values from a uniform sampling $f(x_i)=f( ih), i=1,2,\cdots,n$. This sum approaches the true area, denoted by $S$, of the 
region under the curve $y=f(x)$ in [$a,b]$. The ever improving approximation as $n\to \infty$ is a process of limit and is denoted by 
\be
\lim_{n\to \infty} (b-a)\hat{\bar{f}}[n] = S=I[f,a,b]. 
\ee
Dividing both sides by $b-a$ leads to the limit of $\hat{\bar{f}}[n]$ as $n \to \infty$
\be
\lim_{n\to\infty} \hat{\bar{f}}[n]= \bar{f}.
\ee
As mentioned earlier, the existence of this limit is guaranteed by LLN, and the probability for this limit to fail is zero. 

Therefore, the geometric meaning of the integral $I[f,a,b]$ is the area of the region bounded by 
$y=f(x)$, $y=0$, $x=a$ and $x=b$.

Now following Shen and Lin (2014)[ref. [{\bf 9}]), we can define antiderivative and introduce a concept called derivative-antiderivative (DA) pair. 
We define the integral $I[f,a,x]$ as {\it antiderivative} of $f(x)$ and is denoted by $F(x)$, where $a$ is an arbitrary 
constant and $x$ is  regarded as a variable. We also call 
$f(x)$ the derivative of $F(x)$ and is denoted by $F'(x)=f(x)$. Because $x$ is regarded as a variable, antiderivative  $F(x)$ is a function. 
The function pair $(f(x), F(x))$ is called a {\it DA pair}. 

From the geometric meaning of integral,  $F(x)=I[f,0,x]$ is the area of the region 
bounded by $y=f(t)$, $y=0$, $t=0$ and $t=x$ over the $t-y$ plane (see Figure 3). 

\noindent {\bf Example 1.} Evaluate the antiderivative of $f(x)=1$. 

The area of the rectangle bounded $y=1$, $y=0$, $t=0$ and $t=x$ over the $t-y$ plane is $x$. Thus, $F(x) = x$, and $(1,x)$ is a DA pair.

\noindent {\bf Example 2.} Evaluate the  antiderivative of $f(x)=x$.

The area of the triangle bounded $y=t$, $y=0$, $t=0$ and $t=x$ over the $t-y$ plane is $x^2/2$ (See Figure 2). Thus, $F(x) = x^2/2$, and $(x,x^2/2)$ is a DA pair,
i.e., $(x^2/2)'=x$ and $I[t,0,x]=x^2/2$.

\noindent  {\bf Example 3.} Evaluate the  antiderivative of $f(x)=x^2$. 

This is a problem of calculating the area of a curved triangle (See Figure 3). We can use the uniform sampling as we have done above for  $y=t^2$ in the 
interval $[0,1]$, but now for the interval $[0,x]$. 
\[
\hat{\bar{f}}[n]= \frac{\sum_{i=1}^{n}(xi/n)^2}{n}=x^2\frac{\sum_{i=0}^{n}i^2}{n\times n^2}= x^2\frac{n(n+1)(2 n +1)/6}{n^3}
=x^2 \frac{n+1}{n}\times\frac{2n+1}{6n}.
\]
This mean approaches $x^2/3$ as $n\to \infty$. Hence, 
\[
F(x)=(x-0)\times x^2/3 = x^3/3.
\]
We thus have a DA pair $(x^2, x^3/3)$, i.e., $(x^3/3)'=x^2$ and $I[t^2,0,x]=x^3/3$.

It is tedious to find an antiderivative by definition as shown above. Fortunately, many free software packages, such as WolframAlpha,
 are now available over the internet. 
For example, in the pop up  box from WolframAlpha, 
type  {\tt integral x}\^{\tt 4}.
Press the enter key. The WolframAlpha returns the antiderivative $x^5/5$. 
So $(x^4, x^5/5)$ is a DA pair. 
In the same way, WolframAlpha can generate the following DA pairs.
\begin{enumerate}[(i)]
\item Power function: $(x^n, x^{n+1}/(n+1))$.
\item Exponential function: $(e^x, e^x)$.
\item Natural logarithmic function: $(\frac{1}{x}, \ln x)$.
\item Sine function: $(\cos x, \sin x)$.
\item Cosine function: $(-\sin x, \cos x)$.
\item Tangent function: $(\sec^2 x, \tan x)$.
\end{enumerate}

One can also use {\tt www.wolframalpha.com} to find derivatives. For example, in the pop up box, 
 type {\tt derivative x}\^{\tt 5/5} and press enter. It returns $x^4$. This verifies 
that $(x^5/5)'=x^4$. 

\vskip 0.8cm

\section{Calculation of $I[f,a,b]$} 

\vskip 0.3cm

\noindent We can use antiderivative to calculate an integral (again following the DA pair idea of Shen and Lin (2014) (ref. {\bf 9}]). From the area meaning of an integral (see Figure 4), we have 
the following formula:
\be
I[f,c,d] = F(d)-F(c).
\ee
This is also denoted by $I[f,c,d]=F(x)|^d_c$.
Namely, the area of the region bounded by $y=f(x)$, $y=0$, $x=c$, and $x=d$ is equal to the difference of the area 
bounded by $y=f(x)$, $y=0$, $x=a$, and $x=d$ minus the area bounded by $y=f(x)$, $y=0$, $x=a$, and $x=c$. This is  shown in Figure 4:
$I[f,c,d]$ is the area CDD$'$C$'$, equal to area ODD$'$O$'$minus area OCC$'$O$'$. 

\begin{figure} [ht]
\centering
\includegraphics[height=3in, width=4in]{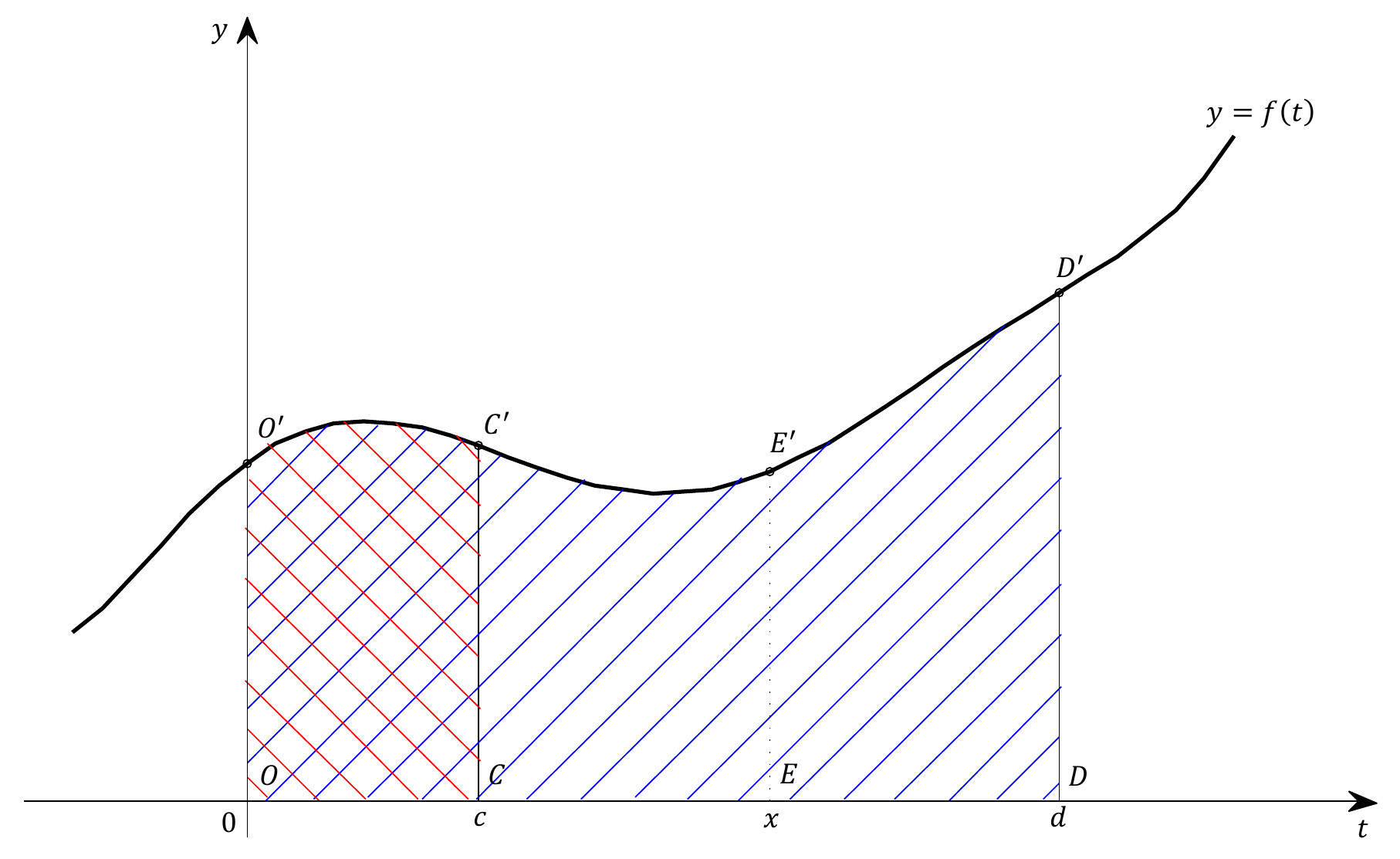}
\caption{Areas and their differences under a curve: an illustration for FTC.}
\label{figure4}
\end{figure}

In traditional calculus, this way of calculating an integral is called Part II of the Fundamental Theorem of Calculus (FTC), and 
the DA pair is called Part I of FTC (e.g., Stewart (2008) [{\bf 10}]).

\noindent {\bf Example 4.} Evaluate $I[x^2,0,1]$.

Since $F(x)=I[x^2] = x^3/3$, $I[x^2,0,1]=x^3/3|_0^1=1^3/3 -0^3/3=1/3.$

\noindent {\bf Example 5.} Evaluate $I[sin^2(x),0,\pi]$.  

In {\tt www.wolframalpha.com}, the command {\tt integrate sin}\^{\tt 2 x} yields  $F(x)=(1/2)(x-\sin x\cos x)$. Thus,
\[
I[\sin^2 x, 0,\pi] = (1/2)(x-\sin x\cos x)|^\pi_0= (1/2)(\pi-\sin \pi\cos \pi)- (1/2)(0-\sin 0\cos 0) = \pi/2.
\]

One can also use  {\tt www.wolframalpha.com} to calculate the integral directly by using the command
{\tt integrate [sin}\^{\tt 2 x, 0, pi]} or {\tt integral [sin}\^{\tt 2 x, 0, pi]}.  This command directly returns value $\pi/2$.

\section{ Use average speed and graphic mean to interpret the meaning of a derivative}

\noindent In the above we have successfully introduced the concepts and calculation techniques of integral and derivative 
from statistics perspective. Further, the meaning of integral of a  function in $[a,b]$
is well interpreted as the value increment of  its anti-derivative function from $x=a$  to $x=b$.  
This section attempts a remaining task of interpreting the meanings of derivative of a function and makes connections 
with the conventional way of defining a derivative. 

 If one drives along a freeway 
from Exit 6 at 2pm and arrives at Exit 98 at 4pm, one's average speed 
is $(98-6)/(4-2)=46$ mph. In general, we use $s(t)$ to represent the location of the car at time $t$,  then the distance traveled 
by the car from time $t_1$ to time $t_2$ is $s(t_2)-s(t_1)$. The average speed is 
\be
\bar{v}= \frac{s(t_2)-s(t_1)}{t_2 - t_1}.
\ee
This kind of average is called a {\it graphic mean}, in contrast to the arithmetic mean discussed earlier. 
Figure 5 gives a schematic illustration of function $y=s(t)$. The graphic mean is thus the slope of the 
secant line, i.e., the purple line that connects $P_1$ and $P_2$ in Figure 5. Namely, $\bar{v}=\tan \theta$. 

\begin{figure} [ht]
\centering
\includegraphics[height=3in,width=4in]{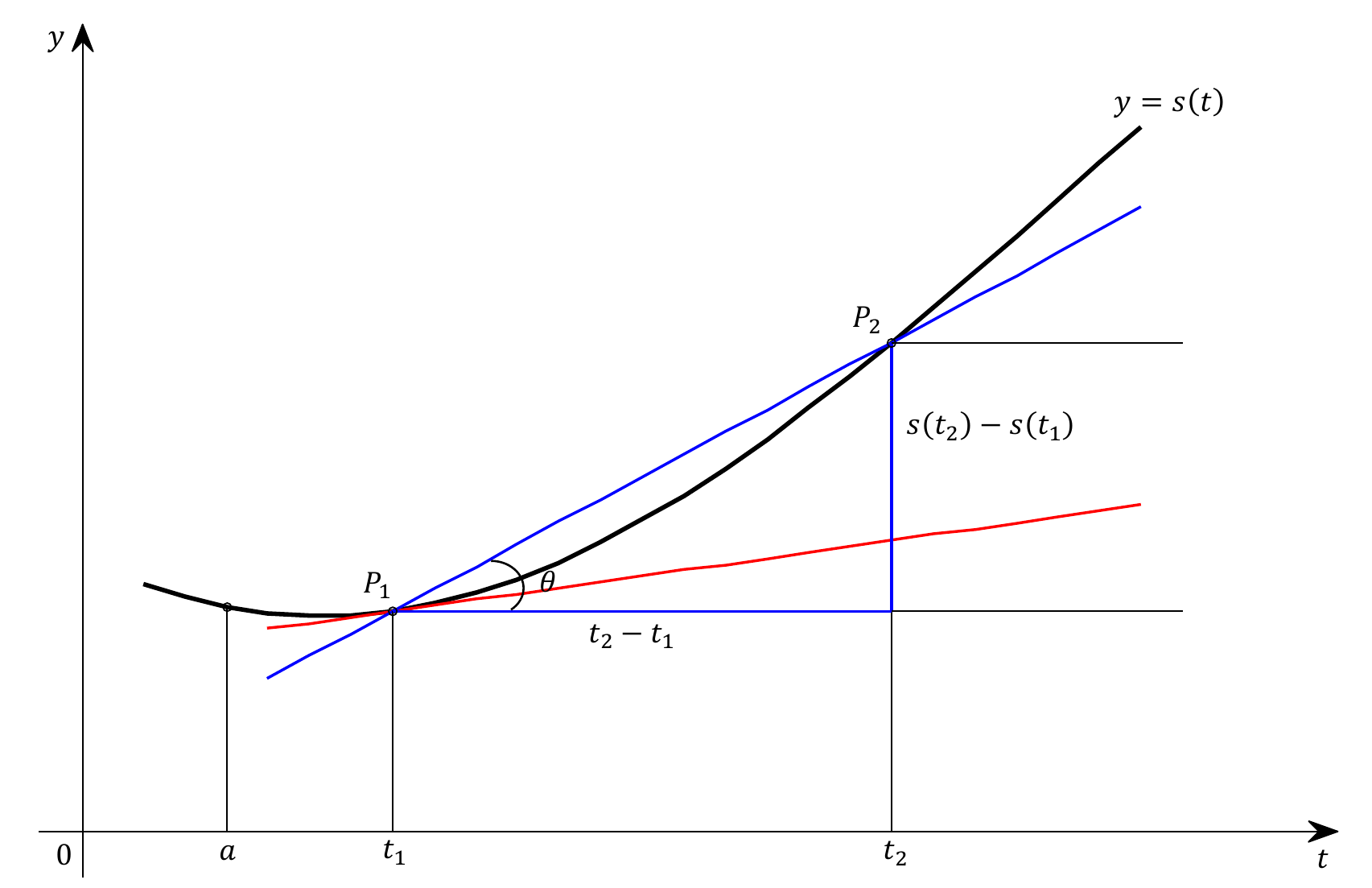}
\caption{Illustration of a graphic mean.}
\label{figure5}
\end{figure}

The instant speed at a time, say $t_1$, is approximately the average speed 
in an ever smaller interval $[t_1, t_2]$, which means $t_2$ goes to $ t_1$. The limit expression is 
\be
\lim_{t_2\to t_1}\frac{s(t_2)-s(t_1)}{t_2 - t_1}= v(t_1).
\label{eq23}
\ee
This is regarded as the {\it slope} of $y=s(t)$ at $t_1$, i.e., the slope of the tangent line (i.e., the red line in Figure 5) 
of $y=s(t)$ at  $t_1$. We would like to show that $s(t)$ is an antiderivative of $v(t)$,
i.e., $s(x)=I[v(t), 0,x] $. We divide $[0,x]$ uniformly into $n$ sub-intervals of width $h=x/n$. Within interval $[x_i, x_{i+1})$, the distance traveled
is approximately $v(x_i)h$. The total distance traveled in $[0,x]$ is $\sum_{i=1}^n v(x_i)h$, which goes to $s(x)-s(0)=I[v,0,x]$ as $n\to\infty$. 
Thus,  $(v(t), s(t))$ is a DA pair, i.e., the derivative of distance with respect to time is speed, and the integral of speed is 
the total distance traveled. Equation (\ref{eq23}) provides another way of defining derivative in addition to the derivative 
introduced via DA pair. This definition of derivative by limit in (\ref{eq23}) is the most popular way of introducing derivatives in today's classrooms.
It is also  the main idea of  Isaac Newton (1642-1727)'s approach to fluxion (see [{\bf 7}] for his book of {\it Method of Fluxion}, 1736).
Before Newton, Pierre de Fermat (1601-1665) developed a  method of tangents in 1629 using a small 
increment (i.e., infinitesimal) which is finally set to be zero when the infinitesimal disappears from the denominator (see [{\bf 2, 3}] ). 
Newton's method of fluxion or tangent is similar to Fermat's but introduced the concept of limit, although he did 
not explicitly use a mathematical notation or formula to express the  limit [See [{\bf 8}] for his book ``{\it The Mathematical 
Principles of Natural Philosophy}" (1729, p45)].

In summary, a derivative of $F(x)$ is a limit of its graphic mean. One can use the formulas of DA 
pairs given in the last section to find the derivative of a given function. For example, to find the derivative of $x^2$, we use 
the DA  pair $(x^n, x^{n+1}/{n+1})$. With $n=1$, the pair becomes $(x, x^2/2)$. So $(x^2/2)'=x$, hence $(x^2)' = 2x$.

\section{ Conclusions and discussion}

 \noindent We have used arithmetic mean to define integral. This definition leads to the DA pair concept, integral's area meaning, 
and FTC. We have used graphic mean to interpret the meanings of a derivative as the slope of a curve or speed of moving 
object. The graphic mean also provides another way of calculating derivative by limit, an application of Fermat's method 
of tangents. Computer calculation for DA pairs is used here in place of traditional derivation of derivative formulas using 
limit and graphic mean. 
 Before today's popularity of computer and internet, the method of introducing derivative by limit 
was very useful in the calculation of calculus problems, since the method can easily be used to derive various kinds of 
derivative formulas for non-polynomial functions, compared to the DA pair approach $I[f,0,x]$. However, 
with today's easy access to notebook computers, smart phones, and  
publicly available software, finding derivatives and integrals for a given function can be readily done 
on the internet. Thus, today's calculus teaching may shift its emphasis from 
the limit-based hand calculations to  the web-based calculations, concepts,  interpretations, and applications. Limit 
 is a difficult, subtle, puzzling and intermediate  procedure for calculus. Although in our statistics 
 calculus here the limit notation is introduced to state the existence of the arithmetic mean 
 following LLN, we hardly used   the concept of limit in formulation derivations and actual calculations.  
 Statistical calculus' rigorous background needs  LLN, which requires a rigorous proof but its conclusion
 has a certain degree of intuition and can be treated as an axiom like those in Euclidean geometry. 
Thus, the calculus based on  LLN implies that, with
today's handy computer technology,
the calculus basics can be axiomized and the calculus method 
can be developed directly with minimum and intuitive theory 
without the need of some unnecessary, intermediate and complex procedures. In this way, calculus can be properly taught 
in high schools or in a workshop of a few hours as proposed in Ref. [{\bf 6}].

\vskip 0.5cm

\noindent  {\bf Acknowledgments.} ~ 
 This research was supported in part by  the US National Science Foundation. The opinions expressed do not necessarily reflect the
views of the funding agencies.
Nancy Tafolla helped calculate the data in Table 2. Discussions with Qun Lin of Chinese Academy of Sciences
were useful for the development of this paper. 
\vskip 1.0cm

\noindent {\Large\bf References}

\begin{enumerate}

\item H. Anton, I.C.  Bivens, and S. Davis, {\it Calculus: Early Transcendentals, Single Variable}, 9th ed. John Wiley and Sons, New York, 880pp, 2008.

\item F. Cajori, {\it A History of Mathematics} (pp. 162-198). 4th ed., Chelsea Publishing Co., New York, 534pp, 1985.

\item D. Ginsburg, B. Groose, J. Taylor, and B. Vernescu,  The History of the Calculus and the Development of Computer Algebra Systems. {\it Worcester Polytechnic 
Institute Junior-Year Project},  www.math.wpi.edu/IQP/BVCalcHist/calctoc.html, 1998

\item R. A. Johnson and G.K. Bhattacharyya,  {\it Statistics: Principles and Methods.} 3rd ed., John Wiley and Sons, New York, 720pp, 1996.

\item D. Kaplan, D. Flath, R. Pruim, and E. Marland,  {\it MAA Ancillary Workshop: Teach Modeling-based Calculus}, at Sheraton Boston on January 3rd, 2012.
\begin{verbatim} http://www.causeweb.org/workshop/jmm12_modeling/. \end{verbatim}

\item Q. Lin, {\it Calculus for High School Students: from a Perspective of Height Increment of a Curve}. People's Education Press, Beijing, 76pp, 2010.

\item I. Newton, {\it The Method of Fluxions and Infinite Series with Application 
to the Geometry of Curve-Lines}. Translated from Latin to English by J. Colson. Printed for Henry Woodfall, London, 339pp, 1736.
\begin{verbatim}
http://books.google.com/books?id=WyQOAAAAQAAJ&printsec=frontcover&source
=gbs_ge_summary_r&cad=0#v=onepage&q&f=false
\end{verbatim}

\item I. Newton,   {\it The Mathematical 
Principles of Natural Philosophy.} Translated from Latin to English by A. Motte. Printed for Benjamin Motte, London, 320pp, 1729.  
\begin{verbatim}
http://books.google.com/books?id=Tm0FAAAAQAAJ&printsec=frontcover&source
=gbs_ge_summary_r&cad=0#v=onepage&q&f=true
\end{verbatim}

\item S.S.P. Shen,  and Q. Lin,  Two hours of simplified teaching materials for direct calculus, {\it Mathematics Teaching and Learning}, No. 2, 2-1 - 2-6, 2014.  arXiv:1404.0070. 

\item J. Stewart,  {\it Single Variable Calculus - Early Transcentals}. 6th ed.,  Thomson Brooks/Cole, Belmont, 763pp, 2008.

\item D.D. Wackerly, W. Mendenhall III, and R.L. Scheaffer,  {\it Mathematical Statistics with Applications.}  6th ed., Duxbury, 853pp, 2002.

\item D. Zazkis, C. Rasmussen, and S.S.P. Shen, A mean-based approach for teaching the concept of integration, {\it PRIMUS (Problems, resources, and Issues in Mathematics Undergraduate Studies)}, {\bf 24}, 116-137, 2014.

\end{enumerate}

\end{document}